\documentclass[11pt, a4paper]{article}
\usepackage[utf8]{inputenc}
\usepackage[T2A]{fontenc}
\usepackage{mathtext}
\usepackage{mathrsfs}
\usepackage{amssymb}
\usepackage{amsfonts}
\usepackage{amscd}
\usepackage{amsmath, amsthm}
\usepackage{mathtools}
\usepackage{hyperref}
\usepackage[symbol]{footmisc}

\DeclareMathOperator*{\res}{Res}

\begin{document}

\begin{center}\textbf{ON POLYNOMIALS OF BINOMIAL TYPE, RAMANUJAN-SOLDNER CONSTANT AND INVERSE LOGARITHMIC DERIVATIVE OPERATOR}
\end{center}
\begin{center}
           \begin{footnotesize}
	DANIL KROTKOV\footnote[2]{Higher School of Economics, Faculty of Mathematics, e-mail: dikrotkov@edu.hse.ru}
	\end{footnotesize}
\end{center}

\begin{footnotesize}
\noindent\textsc{Abstract}. In this work we introduce interesting infinite series, related to Ramanujan-Soldner constant. Our method uses general properties of polynomials of binomial type and Lagrange inversion theorem. Also we study properties of the operator $1/\mathrm{dlog}$, acting on formal power series. In addition, several properties of polynomials of binomial type associated to elementary functions are discussed.
\end{footnotesize}

\renewcommand{\contentsname}{Contents}
\tableofcontents

\begin{footnotesize}
\section*{Introduction}

The function $\mathrm{li}(x)=\int_{0}^{x} dt/\ln t$ is called a \textit{logarithmic integral function}. Many mathematicians studied this function primarily in the field of prime number theorem. But Ramanujan was one of the first, who noticed this function has a unique zero $\mu$ on the real line and calculated the digits of this number in order to check, if it is a known one (see \cite{RamSol}). Since then different mathematicians only computed the new digits without providing any infinite series, converging to this constant, as for example famous Leibniz series for $\pi$: $\frac{\pi}{4}=\sum_{n=0}^\infty \frac{(-1)^n}{2n+1}$.\\
In this work we study the functional inverse to the function $\exp(\mathrm{Ei}(x)-\gamma)$ in order to obtain interesting infinite series including the one, converging to $\mu$, the so-called \textit{Ramanujan-Soldner constant} (here $\gamma$ stands for Euler-Mascheroni constant and $\mathrm{Ei}(x)=\mathrm{li}(e^x)$).\\
When one wants to study $\varphi(x)$ an inverse to an arbitary function $f(x)$, it is often useful to look at polynomials of binomial type generated by this function, i.e. polynomials with property $p_n(\alpha+\beta)=\sum_{k=0}^{n} \binom{n}{k} p_k(\alpha)p_{n-k}(\beta)$ and with generating function $\exp(\alpha\varphi(x))=\sum_{0}^\infty p_n(\alpha)x^n/n!$. One of the useful facts is that an operator $f(\partial/\partial\alpha)$ is a derivative operator for this sequence of polynomials, i.e. $f(\partial/\partial\alpha)p_n(\alpha)=np_{n-1}(\alpha)$, as for example in case of falling factorials $(\alpha)_n=\alpha(\alpha-1)...(\alpha-n+1)$, $(\exp(\partial/\partial\alpha)-1)(\alpha)_n=n(\alpha)_{n-1}$ and $\sum_0^\infty (\alpha)_nx^n/n!=\exp(\alpha\ln(1+x))$.\\
The problem in studying the polynomials, generated by the inverse of $g(x)=\exp(\mathrm{Ei}(x)-\gamma)$ is that the action of operator $g(\partial/\partial\alpha)$ is very difficult to understand compared to the action of $(\exp(\partial/\partial\alpha)-1)f(\alpha)=f(\alpha+1)-f(\alpha)$. In order to avoid such a problem we study the operator $\mathfrak{T}\coloneqq 1/\mathrm{dlog}$ acting on $x+x^2\mathbb{C}[[x]]$ $(f(x) \to f(x)/f'(x))$.\\
In this paper some different topics are discussed. In the first section we present the values of function $M(s)$, considered in the previous paper, at negative integral points and its residues at half-integral points. In the second section we introduce the main properties of inverse of $\exp(\mathrm{Ei}(x)-\gamma)$ and its Taylor coefficients when expanded around zero. In the third section the general construction of polynomials of binomial type is presented and several classical results are disscussed in order to provide the motivation for the questions, raised in section 4. We also look at the basic properties of operator $\mathfrak{T}$ and consider the sequence $...\to f \to f/f' \to ff'/(f'^2-ff'') \to...$ (which we call $\mathfrak{T}$-\textit{chain}). In section 4 we use the methods of section 3 to study functional inverses of specific elementary functions. Also the connection between these functions and the functions introduced in sections 1 and 2 is given. In the final section the equation $\mathfrak{T}^N f(x)=f(x)$ is solved.

\section*{Preliminaries}

\textbf{Key Definitions}
\begin{itemize}
\item Throughout the work the symbol $\gamma$ denotes Euler-Mascheroni constant, i.e. the limit $\lim\limits_{N\to \infty} 1+\frac{1}{2}+$ $+ ...+\frac{1}{N}-\ln N$.
\item $\mathrm{Ei}(x)$ is the exponential integral, $\mathrm{li}(x)$ is the logarithmic integral, $\mu \approx 1.451369...$ is the unique real zero of $\mathrm{li}(x)$:
$$\mathrm{Ei}(x)=-\int_{-x}^{+\infty} \frac{e^{-t}}{t}dt=\ln |x|+\gamma+\int_{0}^{x} \frac{e^t-1}{t}dt$$
$$\mathrm{li}(1+x)=\mathrm{Ei}(\ln(1+x))=\int_{0}^{1+x} \frac{dt}{\ln t}=\ln|x|+\gamma+\int_{0}^{x}\left[ \frac{1}{\ln(1+t)}-\frac{1}{t}\right]dt$$
$$\mathrm{li}(\mu)=\mathrm{Ei}(\ln \mu)=0$$
\item Sequence of polynomials of binomial type is a sequence of polynomials $\{p_n(\alpha)\}_0^\infty$, such that the highest term of $p_n(\alpha)$ is equal to $\alpha^n$ and the convolution property $p_n(\alpha+\beta)=\sum_{k=0}^{n} \binom{n}{k} p_k(\alpha)p_{n-k}(\beta)$ holds (see \cite{Dlt}, \cite{Gl}). As a corollary, their generating function is of the form
$$\exp(\alpha\varphi(x))=\sum_{n=0}^\infty \frac{p_n(\alpha)x^n}{n!}$$
for some Taylor series $\varphi(x) \in x+x^2\mathbb{C}[[x]]$.
\item By $f^{inv}(x)$ we denote the inverse function of $f(x)$, i.e. $f(f^{inv}(x))=f^{inv}(f(x))=x$.
\item The symbol $\buildrel ?\over =$ is used when the convergence of the series is unknown, but in case of convergence it has the suggested limit.
\end{itemize}
\textbf{Key Tools}
\begin{itemize}
\item \textit{Lagrange inversion theorem} (see \cite{Lag}). Consider the function $f(x)=\sum_{n=1}^\infty a_n x^n$ and $a_1=1$. Then Taylor series of its inverse is of the form
$$f^{inv}(x)=\sum_{n=1}^\infty \frac{x^n}{n}\res_{t=0} \frac{1}{f^n(t)}dt$$
\item For any power series $f(x)=\sum_{n=0}^\infty a_n x^n$ the action of operator $f(\partial/\partial\alpha)$ is well-defined on formal series of the form $\sum_{n=0}^\infty b_n \alpha^{s-n}$ for any $s$ and well-defined on exponents $e^{\alpha C}$, whenever $f(x)$ converges at point $C$, so $f(\partial/\partial\alpha)e^{\alpha C}=f(C)e^{\alpha C}$. As a consequence, its action is well-defined on polynomials. For polynomials of binomial type generated by the function $\varphi(x)=f^{inv}(x)$, $e^{\alpha \varphi(x)}=\sum_{n=0}^\infty p_n(\alpha)x^n/n!$, there is an obvious shift property $f(\partial/\partial\alpha)p_n(\alpha)=np_{n-1}(\alpha)$, which can be verified directly.
\item Suppose $\varphi(x)=f^{inv}(x)$ generates polynomials $p_n(\alpha)$. Consider $\gamma(x)=(f(x)e^{-x})^{inv}$. Then according to Lagrange inversion theorem
$$e^{\alpha \gamma(x)}=\sum_{n=0}^\infty \frac{\alpha  p_n(\alpha +n)}{\alpha +n}\frac {x^n}{n!}$$
\end{itemize}
One of the main goals of this paper was to understand, which deformations of the initial function $f(x)$ in some sense preserve power series of the inverse functions. The transform $f(x) \to f(x)e^{-x}$ has been known at least since Abel (see \cite{Abl}), and as we see, acts nicely on polynomials generated by their inverses. But to understand the motivation behind the study of transform $f(x)\to f(x)/f'(x)$ we have to look at this useful lemma.
\begin{itemize}
\item Proposition ($\mathfrak{T}$-\textit{operator}):
$$\frac{f}{f'}\left(\frac{\partial}{\partial \alpha}\right) p_n(\alpha)=n \frac{p_{n}(\alpha)}{\alpha}$$
\textit{Proof}:
\begin{align*}
\sum_{n=0}^\infty \frac{x^n}{n!}\frac{f}{f'}\left(\frac{\partial}{\partial \alpha}\right) p_n(\alpha)&=\frac{f}{f'}\left(\frac{\partial}{\partial \alpha}\right) e^{\alpha \varphi(x)}=x\varphi'(x)e^{\alpha \varphi(x)}=\\
\tag*{\qed}&=\frac{x}{\alpha} \frac{\partial}{\partial x}e^{\alpha \varphi(x)}=\sum_{n=1}^\infty \frac{p_n(\alpha)}{\alpha} \frac{x^n}{(n-1)!}
\end{align*}
\end{itemize}
This formula shows that even when we don't understand how the operator $f(\frac{\partial}{\partial\alpha})$ acts, but understand the behaviour of the functions under the action of $f(\frac{\partial}{\partial\alpha})/f'(\frac{\partial}{\partial\alpha})$, the polynomials of such kind are expected to have interesting properties.
\end{footnotesize}

\newpage

\hrule

\section{Integral representation for $M(s)$ and values of its analytic continuation at negative integral points}

\hrule

~\\
~\\
\underline{\textbf{Theorem 1.1}}: Consider the function
$$
M(s)=\sum\limits_{n=1}^\infty \dfrac{e^{-n}}{n^s}\left( 1+n+\dfrac{n^2}{2!}+...+\dfrac{n^n}{n!} \right)
$$
Then it has the representation
$$
\frac{1}{2}\Gamma(1+s)M(s)=\int_{0}^{1} \frac{ (-t-\ln(1-t))^s}{t^3}dt
$$
The proof is divided in a few parts. Consider the Lambert $W$-function which has the properties:
$$
W(x)=(xe^x)^{inv} \Rightarrow e^{\alpha W(x)}=\sum_{n=0}^\infty \frac{\alpha(\alpha-n)^{n-1}x^n}{n!};~~ W'(x)=\frac{W(x)}{x(1+W(x))}
$$
\textbf{Proposition 1.1}:
$$
\sum\limits_{k=0}^{n} \frac{n^k}{k!}=\int\limits_{0}^{+\infty} \frac{(n+t)^n}{n!}e^{-t}dt\eqno (1.1)
$$
\emph{Proof}:
\begin{align*}
\sum\limits_{k=0}^{n} \frac{n^k}{k!}=\sum_{k=0}^{n}\binom{n}{k}\frac{n^k}{n!}(n-k)!&=\sum_{k=0}^{n}\binom{n}{k}\frac{n^k}{n!}\int\limits_{0}^{+\infty} t^{n-k}e^{-t}dt\\
\tag*{\qed}&=\int\limits_{0}^{+\infty} \frac{(n+t)^n}{n!}e^{-t}dt
\end{align*}
\textbf{Proposition 1.2}:
$$
\sum_{n=0}^\infty x^n \sum\limits_{k=0}^{n} \frac{n^k}{k!}=\frac{1}{(1+W(-x))^2}\eqno (1.2)
$$
\emph{Proof}:
\begin{align*}
\sum_{n=0}^\infty x^n\sum\limits_{k=0}^{n} \frac{n^k}{k!}&=\int\limits_{0}^{+\infty} e^{-t} \sum_{n=0}^\infty \frac{x^n (t+n)^n}{n!}dt =&\\
&=\int\limits_{0}^{+\infty} e^{-t} \left( \sum_{n=0}^\infty \frac{t(t+n)^{n-1} x^n}{n!}+ \sum_{n=0}^\infty \frac{n(t+n)^{n-1} x^n}{n!}\right)dt=&\\
&=\int\limits_{0}^{+\infty} e^{-t} \left( e^{-tW(-x)}+\frac{x}{t}\frac{d}{dx}e^{-tW(-x)}\right)dt=
\end{align*}
\begin{align*}
&=(1+xW'(-x))\int\limits_{0}^{+\infty} e^{-t(1+W(-x))}dt=\\
\tag*{\qed}&=\frac{(1+xW'(-x))}{1+W(-x)}=\frac{1}{(1+W(-x))^2}
\end{align*}
\textbf{Proposition 1.3}:
$$
\sum_{n=1}^\infty \frac{x^n}{n} \sum\limits_{k=0}^{n} \frac{n^k}{k!}=-W(-x)-\ln(1+W(-x))\eqno (1.3)
$$
\emph{Proof}: Take the derivative of (1.3) and use the differential equation of $W(x)$, to obtain (1.2). The result follows after comparing the values of LHS and RHS of (1.3) at $0$. \qed\\
~\\
\textbf{Remark 1.1}. The latter formula allows us to evaluate the following limit in much more easier manner, than it was demonstrated in the previous paper:
\begin{align*}
\sum_{n=1}^\infty \left[ \frac{e^{-n}}{n}\sum_{k=0}^{n}\frac{n^k}{k!}-\frac{1}{2n}\right]&=\lim_{x \to 1^{-}} -W\left(-\frac{x}{e}\right)-\ln\left(1+W\left(-\frac{x}{e}\right)\right)+\frac{1}{2}\ln(1-x)\\
&=1-\frac{1}{2}\ln 2
\end{align*}
\textbf{Remark 1.2}. It follows from (1.2) that $M(s)$ has the following integral representation in terms of Mellin transform:
$$
\Gamma(s)M(s)=\int\limits_{0}^{+\infty} t^{s-1}\left( \frac{1}{(1+W(-e^{-1-t}))^2} -1\right) dt\eqno (1.4)
$$
Now it is easy to derive from (1.4) the statement of Theorem 1.1.\\

\noindent\textbf{Proposition 1.4}:
\begin{align*}
\tag{1.5}\Gamma(s)M(s)&=\int_{0}^{1} (-t-\ln(1-t))^{s-1} \left(1+\frac{1}{t}\right) dt\\
\tag{1.6}&=\frac{2}{s}\int_{0}^{1} \frac{ (-t-\ln(1-t))^s}{t^3}dt
\end{align*}
\emph{Proof}: The first formula follows from (1.4) after the change of variable $t \rightarrow -t-\ln(1-t)$, since $W(xe^x)=x$ on the integration domain. The second equality follows from the first one via integration by parts. The claim of the theorem follows immediately.\hfill\ensuremath{\blacksquare}\\
~\\
This integral representation is useful, because it allows one to derive the values of $M(-k)$ in terms of the values of special polynomials. In order to do that consider the following generating function:
\begin{align*}
& \left[ \frac{-t-\ln(1-t)}{t^2/2}\right]^s=\left[ 1+\frac{2t}{3}+\frac{2t^2}{4}+\frac{2t^3}{5}+...\right]^s=\sum_{k=0}^\infty A_k(s)t^k\\
&\left( A_0(s)=1, ~A_1(s)=\frac{2s}{3}, ~A_2(s)=\frac{s(4s+5)}{18},...\right)
\end{align*}
\underline{\textbf{Theorem 1.2}}: Suppose $N \in \mathbb{Z}_{\geqslant 0}$. Then the following holds true:
\begin{align*}
&0) M(0)=A_2'(0)-A_0(0)=-\frac{13}{18};\\
&1) M(-N)=(-1)^{N-1} 2^N (N-1)!A_{2N+2}(-N), ~\text{if} ~(N>0);\\
&2) \lim_{s \to -N} M\left(s+\frac{1}{2}\right)(s+N)=(-1)^{N-1} \frac{(2N)!}{2^N N!} \frac{A_{2N+1}(\tfrac{1}{2}-N)}{2N-1}\sqrt{\frac{2}{\pi}}.
\end{align*}
\emph{Proof}: One can rewrite the formula (1.6) in the following way:
\begin{align*}
2^{s-1}\Gamma(s+1)M(s)&=\int_{0}^{1} t^{2s-3}\left[ \frac{-t-\ln(1-t)}{t^2/2}\right]^s dt\\
\tag{1.7}&=\sum_{k=0}^\infty A_k(s)\int_{0}^{1}t^{2s-3+k}dt=\sum_{k=0}^\infty \frac{A_k(s)}{2s-2+k}
\end{align*}
0) Substituting $s=0$ in both sides of this equality, it is easy to see that only two terms the 0-th one and the 2-nd one would appear to be nonzero, i.e.
$$
\frac{1}{2}M(0)=\frac{A_0(0)}{-2}+\frac{A_2(s)}{2s} \Biggr|_{s=0}=-\frac{1}{2}+\frac{4s+5}{36}\Biggr|_{s=0}=-\frac{13}{36}
$$
(note: here we used the explicit formula for polynomial $A_2$ which was presented on the previous page) \qed\\
~\\
1) Note that whether the series (1.7) converges or not, the following holds true:
\begin{align*}
\lim_{s \to -N} 2^{s-1} \frac{\Gamma(s+1)M(s)}{\Gamma(2s-2)}&=\lim_{s \to -N} \sum_{k=0}^\infty \frac{A_k(s)(2s-2)(2s-1)...(2s+k-3)}{\Gamma(2s+k-1)}\\
~\\
&=\lim_{s \to -N} A_{2N+2}(s)(2s-2)(2s-1)...(2s+2N-1)
\end{align*}
This equality holds because the first $2N+1$ terms are equal to zero, since $\Gamma$-function has poles at negative integer points, and the terms going after the one with index $2N+2$ vanish because of the zeroes of falling factorials. It is left to say that the following limit exists:
$$
\lim_{s \to -N} \frac{\Gamma(s+1)}{\Gamma(2s-2)}
$$
to obtain the desired result. \qed\\
~\\
2) Similarly,
$$
\lim_{s \to \tfrac{1}{2}-N} 2^{s-1} \frac{\Gamma(s+1)M(s)}{\Gamma(2s-2)}=-A_{2N+1}\left(\frac{1}{2}-N\right) (2N+1)!
$$
It is left to say that the limit
$$
\lim_{s \to \tfrac{1}{2}-N} \frac{\Gamma(s+1)}{\Gamma(2s-2)(s+N-\tfrac{1}{2})}
$$
exists to obtain the desired result. \hfill\ensuremath{\blacksquare}\\
~\\
\textbf{Remark}. It would be interesting to know if $M(s)$ has poles at all negative half-integral points or not, and does it have zeroes at negative integers or not.

\newpage

\hrule

\section{A few observations on Ramanujan-Soldner constant and related number pyramid}

\hrule

~\\
~\\
\underline{\textbf{Theorem 2.1}}: Consider Ramanujan-Soldner constant $\mu$, $\mathrm{Ei}(\ln \mu)=0$, and $\mathrm{Ei}$ --- exponential integral function. Consider the sequence $a_n$:
$$
a_n\coloneqq\frac{1}{n} \res_{z=0} ~ e^{-n(\mathrm{Ei}(z)-\gamma)}~ dz=\frac{1}{n} \res_{z=0} ~\exp\left(-n\int_{0}^{z} \frac{e^t-1}{t}dt \right) ~ \frac{dz}{z^n}
$$
Then the following holds:
\begin{align*}
&\sum_{n=1}^\infty a_n e^{-\gamma n}\buildrel ?\over =\ln \mu ~&&\sum_{n=1}^\infty (-1)^{n-1}a_n e^{-\gamma n}=\infty\\
&\sum_{n=1}^\infty \frac{a_n}{n} e^{-\gamma n}=\mu -1 ~&&\sum_{n=1}^\infty \frac{(-1)^{n-1}a_n}{n} e^{-\gamma n}=1\\
& ~&&\sum_{n=1}^\infty \frac{(-1)^{n-1}a_n}{n^2} e^{-\gamma n}=\ln 2
\end{align*}
The proof is divided in a few parts. Suppose
$$
\psi(x)\coloneqq\left( x\exp \left( \int_{0}^{x} \frac{e^t-1}{t} dt \right) \right) ^{inv}
$$
Notice that $\psi(z)$ is an inverse of holomorphic function, whose Taylor series around zero is of the form $x+x^2\mathbb{C}[[x]]$. That means $\psi(x)$ also has the expansion of such a form around zero by the Lagrange inversion theorem with some radius of convergence.\\
~\\
\textbf{Proposition 2.1}:
$$
\psi(x)=\sum_{n=1}^\infty a_n x^n \eqno (2.1)
$$
\emph{Proof}: The result obviously follows from the definition of $a_n$ and Lagrange inversion formula. \qed\\
~\\
This approach allows us to obtain a representation of $a_n$, but unfortunately it is too complicated to show some important properties of these numbers.\\
~\\
\textbf{Proposition 2.2} (exact formula):
\begin{align*}
a_n&=\frac{1}{n}\sum_{k=0}^{n-1} \frac{(-n)^k}{k!}\sum_{\substack{m_1+...+m_k=n-1\\ m_i>0}} \frac{1}{m_1 m_1!...m_k m_k!}\\
\tag*{(2.2)}&=\frac{1}{n!}\sum_{k=0}^{n-1} \frac{(-n)^k}{k!}\sum_{\substack{m_1+...+m_k=n-1\\ m_i>0}} \binom{n-1}{m_1,...,m_k}\frac{1}{m_1...m_k}
\end{align*}
\emph{Proof}:
\begin{align*}
a_n&=\frac{1}{n} \res_{z=0} ~\exp\left(-n\int_{0}^{z} \frac{e^t-1}{t}dt \right) ~ \frac{dz}{z^n}\\
&=\frac{1}{n} \res_{z=0} ~\sum_{k=0}^{n-1} \frac{(-n)^k}{k!} \left( \sum_{m=1}^\infty \frac{z^m}{m\cdot m!} \right)^k ~ \frac{dz}{z^n}\\
\tag*{\qed}&=\frac{1}{n}\sum_{k=0}^{n-1} \frac{(-n)^k}{k!}\sum_{\substack{m_1+...+m_k=n-1\\ m_i>0}} \frac{1}{m_1 m_1!...m_k m_k!}
\end{align*}
\textbf{Proposition 2.3} (inversion of the logarithmic integral):
$$
\left(xe^{\mathrm{li}(1+x)-\ln|x|-\gamma} \right)^{inv}=e^{\psi(z)}-1=\sum_{n=1}^\infty \frac{a_n}{n}z^n\eqno (2.3)
$$
\emph{Proof}: Suppose $T(x)=x\exp \left( \mathrm{Ei}(x)-\ln|x|-\gamma \right)=\psi^{inv}$. Then it follows that
\begin{align*}
\tag*{\qed}\frac{T(x)}{T'(x)}=xe^{-x} \Rightarrow x \psi'(x)=\psi(x)e^{-\psi(x)} \Rightarrow (e^{\psi(x)}-1)'=\frac{\psi(x)}{x}
\end{align*}
\textbf{Remark 2.1}. It is interesting how properties of this function allow us to express the coefficients of $(\exp(\mathrm{li}(1+x)-\gamma))^{inv}$ in terms of $(\exp(\mathrm{Ei}(x)-\gamma))^{inv}$ without using the direct exponentiation. Such an approach also allows us to write down the following recurrence relation for $a_n$.\\
~\\
\textbf{Proposition 2.4}:
$$
(1-n)~a_n=\sum_{k=1}^{n-1} \frac{n-k}{k}~a_k a_{n-k}\eqno (2.4)
$$
\emph{Proof}: This recurrence relation corresponds to the differential equation:
$$
x\psi'(x) \int_{0}^{x} \frac{\psi(t)}{t}dt=-x^2 \left( \frac{\psi(x)}{x}\right)'
$$
By (2.3) the following holds:
\begin{align*}
x\psi'(x) \int_{0}^{x} \frac{\psi(t)}{t}dt&=x\psi'(x) (e^{\psi(x)}-1)\phantom{\int_{0}^{x}}\\
&=x\psi'(x)e^{\psi(x)}-x\psi'(x)\phantom{\int_{0}^{x}}\\
&=\psi(x)-x\psi'(x)\phantom{\int_{0}^{x}}\\
\tag*{\qed}&=-x^2 \left( \frac{\psi(x)}{x}\right)'\phantom{\int_{0}^{x}}
\end{align*}
This relation is invariant under the transform $a_n \rightarrow a_n C^n$. So using the property $a_1=1$ and the induction it is easy to show that $(-1)^{n-1}a_n>0$:
$$
(n-1)~(-1)^{n-1}a_n=\sum_{k=1}^{n-1} \frac{n-k}{k}~(-1)^{k-1}a_k ~ (-1)^{n-k-1}a_{n-k}>0
$$
To move on it is necessary to find the radius of convergence of the series $\psi(x)$.\\
~\\
\textbf{Proposition 2.5}:
$$
\forall z, ~|z|<e^{-\gamma}: ~~ \sum_{n=1}^\infty a_n z^n ~ \text{converges}.\eqno (2.5)
$$
\emph{Proof}: $(\psi^{inv})'=T'=\exp(x+\mathrm{Ei}(x)-\ln|x|-\gamma)$. This function has zero only at the point $x=-\infty$. Since $\mathrm{Ei}(-\infty)=0$ we obtain that $T(-\infty)=-e^{-\gamma}$.  This means that Taylor series for the function $\psi(z)$ converges in the region $|z|<e^{-\gamma}$. \qed\\
~\\
\textbf{Remark 2.2}. Suppose $b_n=(-1)^{n-1}a_n e^{-\gamma n}>0$. It follows from the property $T(-\infty)=$ $=-e^{-\gamma}$ that the series $\sum_{1}^{\infty} b_n$ diverges. But it is not that trivial in case of the alternating series. The convergence of the series $\sum_{1}^\infty (-1)^{n-1} b_n$ would follow from the fact that the sequence $b_n$ is strictly decreasing. Unfortunately, the author does't know the proof of this fact (hypothesis 2.1). Also the convergence would follow from the estimate $b_n = O (\frac{1}{n})$ by Littlewood's Tauberian theorem (hypothesis 2.2).\\
~\\
\textbf{Hypothesis 2.1}: $\forall n \in \mathbb{N}: b_n > b_{n+1}$.\\
~\\
\textbf{Proposition 2.6}: $b_n \xrightarrow{n \to \infty} 0$.  \hfill\ensuremath{(2.6)}\\
\emph{Proof}: see Appendix A.\\
~\\
\textbf{Hypothesis 2.2}: The sequence $\{nb_n\}_1^\infty$ is strictly increasing to 1.\\
(note: from (2.6) it only follows that if the limit $\lim_\infty nb_n$ exists, then it is equal to 1)\\
~\\
\textbf{Proposition 2.7}: If the series $\sum_1^\infty (-1)^{n-1} b_n$ converges, then its value is $\ln \mu$.\\
\emph{Proof}: If the series converges, then by Abel theorem the limit is the value of the function at that point. $T(\ln \mu)=e^{-\gamma} \Rightarrow \psi(e^{-\gamma})=\ln \mu$. \qed\\
~\\
\textbf{Remark 2.3}. It follows from (2.6) that $\sum_{n=1}^\infty b_n x^n n^{-\varepsilon}$ converges if $\varepsilon > 1$ for all $x, ~ |x|\leqslant 1$. It means that to finish the proof of the theorem it is left to find the values of the four other introduced infinite series.\\
~\\
\textbf{Proposition 2.8}:
\begin{align*}
\tag{2.7}\Gamma(s+1)\sum_{n=1}^\infty \frac{b_n}{n^s}&=\int_{0}^{+\infty} (-\mathrm{Ei}(-x))^sdx\\
\tag{2.8}\Gamma(s+1)\sum_{n=1}^\infty \frac{(-1)^{n-1}b_n}{n^s}&=\int_{0}^{\ln \mu} (-\mathrm{Ei}(x))^sdx
\end{align*}
\emph{Proof}: The standard representations of these series in terms of Mellin transform follow from the expansion $x e^{-\gamma}\psi'(-e^{-\gamma}x)=\sum_1^\infty nb_n x^n$. After the substitution $t \rightarrow -\ln(T(t)e^\gamma)$ they could be transformed to the desired ones.\\
Now it is left to say that from (2.7-8) one can easily derive the exact value of each of the remaining series using integration by parts.\hfill\ensuremath{\blacksquare}\\
~\\
\textbf{Remark 2.4}. We should also mention the following nice series:
\begin{center}
\boxed{\sum_{n=1}^\infty \frac{b_n}{n^3}=\frac{\pi^2}{6}-\sum_{n=0}^\infty \frac{1}{2^{2n} (2n+1)^2}}
\end{center}
\underline{\textbf{Theorem 2.2}}:
$$
\left[ \frac{d^n}{dx^n} e^{\alpha \mathrm{Ei}(x)}\right] e^{-\alpha \mathrm{Ei}(x)}=\sum_{k=1}^{n} \sum_{m=k}^{n} \frac{(-1)^{k-m} \alpha^k e^{kx}}{x^m} A_{k, m}^n
$$
and the numbers $A_{k, m}^n$ satisfy the relation
$$
A_{k, m}^{n+1}=k A_{k, m}^n+(m-1)A_{k, m-1}^n+A_{k-1, m-1}^n
$$
In particular,
$$
A_{m, m}^n = \begin{Bmatrix} n\\m \end{Bmatrix}; A_{m, n}^n =  \begin{bmatrix} n\\m \end{bmatrix}; A_{1, m}^n = \frac{(n-1)!}{(n-m)!},
$$
where $\begin{Bmatrix} n\\m \end{Bmatrix}$ stands for the absolute values of Stirling numbers of the second kind,\\
and $\begin{bmatrix} n\\m \end{bmatrix}$ stands for the absolute values of Stirling numbers of the first kind.\\
~\\
~\\
\emph{Proof}: For $n=1$:
$$
\left[ \frac{d}{dx} e^{\alpha \mathrm{Ei}(x)}\right] e^{-\alpha \mathrm{Ei}(x)}=\frac{\alpha e^x}{x}
$$
Now use induction on $n$ to obtain
\begin{align*}
\frac{d}{dx} &\left[ e^{\alpha \mathrm{Ei}(x)}\sum_{k=1}^{n} \sum_{m=k}^{n} \frac{(-1)^{k-m} \alpha^k e^{kx}}{x^m} A_{k, m}^n\right]=\\
&= e^{\alpha \mathrm{Ei}(x)}\sum_{k=1}^{n} \sum_{m=k}^{n} \frac{(-1)^{k-m} \alpha^{k+1} e^{(k+1)x}}{x^{m+1}} A_{k, m}^n+\\
&+ e^{\alpha \mathrm{Ei}(x)}\sum_{k=1}^{n} \sum_{m=k}^{n} \frac{(-1)^{k-m} \alpha^k k e^{kx}}{x^m} A_{k, m}^n+\\
&+ e^{\alpha \mathrm{Ei}(x)}\sum_{k=1}^{n} \sum_{m=k}^{n} \frac{(-1)^{k-m-1} \alpha^k m e^{kx}}{x^{m+1}} A_{k, m}^n=\\
&= e^{\alpha \mathrm{Ei}(x)}\sum_{k=1}^{n+1} \sum_{m=k}^{n+1} \frac{(-1)^{k-m} \alpha^k e^{kx}}{x^m} \left(k A_{k, m}^n+(m-1)A_{k, m-1}^n+A_{k-1, m-1}^n \right)
\end{align*}
Setting $k=m$, one can obtain the relation $A_{m, m}^{n+1}=m A_{m, m}^n+A_{m-1, m-1}^n$. Since $A_{1,1}^1=1$, we indeed deal with Stirling numbers of the second kind. The remaining claims can be proved by the same argument.\hfill\ensuremath{\blacksquare}\\
~\\
So the number pyramid, containing $A_{k,m}^n$, has the property that each of its faces coincides with one of fundamental triangular arrays: Pascal's and Stirling triangular arrays.\\
~\\
The slices of this pyramid for $n=1,...,6$ are listed in Appendix B.

\newpage

\hrule

\section{Polynomials of binomial type and Ramanujan's Master Theorem}

\hrule

~\\
~\\
\textbf{Remark}. This section mostly consists of formal and trivial equalities. Polynomials of binomial type are being studied in many articles on combinatorics and umbral calculus, so the major part of conclusions in this section are not original and are given here only to provide the motivation for this research.
\begin{center}
	\textbf{Trivial theorems}
\end{center}
Consider $f(x) \in x+x^2 \mathbb{C}[[x]]$. Suppose $\varphi(x)=f^{inv}(x)$ and $e^{\alpha \varphi(x)}=\sum_{0}^\infty p_n(\alpha)\frac{x^n}{n!}$. Then the convolution formula holds $p_n(\alpha + \beta)=\sum_{k=0}^{n} \binom{n}{k} p_k(\alpha) p_{n-k}(\beta)$. As it was already mentioned in preliminaries, operator $f(\partial/\partial \alpha)$ has a well-defined action on the ring of polynomials and more generally on the formal asymptotic series (series of the form $\alpha^s\mathbb{C}[[\alpha^{-1}]]$). Also $f(\partial/\partial \alpha)e^{C\alpha}=f(C)e^{C\alpha}$.\\ 
~\\
\textbf{Remark 3.1}. Delta operator (see ~\cite{Dlt}):
$$
f\left(\frac{\partial}{\partial \alpha}\right) p_n(\alpha)=n p_{n-1}(\alpha) \eqno (3.1)
$$
\textbf{Remark 3.2}. $\mathfrak{T}$-operator (see preliminaries):
$$
\frac{f}{f'}\left(\frac{\partial}{\partial \alpha}\right) p_n(\alpha)=n \frac{p_{n}(\alpha)}{\alpha} \eqno (3.2)
$$
\textbf{Example}: Consider $\psi(x)$, the function from section 2, i.e.
$$
\psi(x)\coloneqq\left( x\exp \left( \int_{0}^{x} \frac{e^t-1}{t} dt \right) \right) ^{inv}
$$
Consider the expansions
\begin{align*}
&e^{\alpha \psi(x)}=\sum_{n=0}^\infty \frac{\nu_n(\alpha)x^n}{n!}; && e^{\alpha (e^{\psi(x)}-1)}=\sum_{n=0}^\infty \frac{\mu_n(\alpha)x^n}{n!}\\
\text{Then}~~~ & \nu_n '(\alpha-1)=n\frac{\nu_n (\alpha)}{\alpha}; && \int_{0}^\infty \frac{e^{-t}}{t} \left( \mu_n(\alpha)-\mu_n (\alpha-t)\right) dt=n\frac{\mu_n(\alpha)}{\alpha}
\end{align*}
\textbf{Remark 3.3}. Deformation by $e^{-x}$ (see preliminaries): Consider $\gamma(x)=(f(x)e^{-x})^{inv}$. Then
$$
e^{\alpha \gamma(x)}=\sum_{n=0}^\infty \frac{\alpha  p_n(\alpha +n)}{\alpha +n}\frac {x^n}{n!} \eqno (3.3)
$$
Suppose one wants to continue the function $f(n, \alpha)=p_n(\alpha)$ as a function of $n$ to the fractional argument in the same manner as for $f=\varphi=\mathrm{Id}$:
$$
e^{\alpha \varphi(x)}=\sum_{n=0}^\infty \frac{\alpha^n x^n}{n!} ~~\text{and}~~ \alpha^{-s}=\frac{1}{\alpha^s}
$$
or as in case $f(x)=e^x-1$, $\varphi(x)=\ln(1+x)$:
$$
e^{\alpha \varphi(x)}=\sum_{n=0}^\infty \frac{x^n}{n!}(\alpha)_{n} ~~\text{and}~~ (\alpha)_{-s}=\frac{\Gamma(1+\alpha)}{\Gamma(1+\alpha+s)}
$$
In order to do that we can use Ramanujan's Master Theorem (see~\cite{RamMas}):
$$
\Gamma(s)p_{-s}(\alpha)=\int_{0}^{?} t^{s-1}e^{\alpha \varphi(-t)}dt
$$
Changing variables and integrating by parts one could then derive:
$$
\Gamma(s)\frac{p_{1-s}(\alpha)}{\alpha}=\int_{0}^{?} (-f(-t))^{s-1}e^{-\alpha t}dt
$$
It should be mentioned that setting $?=+\infty$ in the first formula,  we obtain incorrect formula for $(\alpha)_{-s}$:
$$
\Gamma(s)(\alpha)_{-s}=\int_{0}^{1} t^{s-1} (1-t)^{\alpha}dt \neq \int_{0}^{+\infty} t^{s-1}e^{\alpha \ln(1+(-t))}dt.
$$
But setting $?=+\infty$ in the second formula, we obtain correct equalities in both cases whenever $\alpha>0, \Re s>0$:\\
\begin{align*}
\Gamma(s)\frac{\alpha^{1-s}}{\alpha}&=\int\limits_{0}^{+\infty} (-(-t))^{s-1}e^{-\alpha t}dt\\
\Gamma(s)\frac{(\alpha)_{1-s}}{\alpha}&=\int\limits_{0}^{+\infty} (1-e^{-t})^{s-1}e^{-\alpha t} dt
\end{align*}
~\\
\textbf{Proposition 3.1}. (\emph{formal}). Suppose $f(x) \in x+x^2 \mathbb{R}[[x]]$ and the series converges everywhere. Suppose $\int_{0}^{+\infty} (-f(-t))^{s-1}e^{-\alpha t}dt$ converges for all $\alpha>0, s>0$. Then the canonical continuation of $p_n(\alpha)$ as a function of $n$ is given by the formula
$$
p_{1-s}(\alpha)\coloneqq\frac{\alpha}{\Gamma(s)}\int\limits_{0}^{+\infty} (-f(-t))^{s-1}e^{-\alpha t}dt \eqno (3.4)
$$
~\\
\textbf{Remark 3.4}. Of course, a continuation of the function $g: \mathbb{N} \to \mathbb{C}$ is not unique, that's why this proposition should be considered as a definition of ``polynomials of binomial type canonical continuation'', not as a theorem.\\
~\\
\textbf{Remark 3.5}.
\begin{align*}
p_{-1}(\alpha)=\alpha\int\limits_{0}^{+\infty} (-f(-t))e^{-\alpha t}dt&=-\alpha f\left(\frac{\partial}{\partial\alpha}\right)\int\limits_{0}^{+\infty}e^{-\alpha t}dt\\
&=-\alpha f\left(\frac{\partial}{\partial\alpha}\right)\frac{1}{\alpha}
\end{align*}
So $p_{-1}(\alpha)$ corresponds to Laplace transform of the function $f$. That means even if the integrals diverge, binomial ``polynomials'' of negative index can be evaluated in the ring of formal Laurent series.\\
~\\
\textbf{Proposition 3.5}. (\emph{formal}).
$$
p_{1-s}(\alpha)=\alpha \left( \frac{f\left(\frac{\partial}{\partial \alpha}\right)}{\frac{\partial}{\partial \alpha}}\right)^{s-1} \alpha^{-s} \eqno (3.5)
$$
\emph{Proof}: The derivation is similar to the previous one. \qed\\
~\\
\textbf{Proposition 3.6}. (\emph{formal}). Suppose $f(x) \in x+x^2\mathbb{R}[[x]]$, $\varphi=f^{inv},\\ \gamma(x)=(f(x)e^{-x})^{inv}$. Suppose
$$
e^{\alpha \varphi(x)}=\sum_{n=0}^\infty \frac{p_n(\alpha)x^n}{n!};~~~~~ e^{\alpha \gamma(x)}=\sum_{n=0}^\infty \frac{q_n(\alpha)x^n}{n!}
$$
Then the following holds:
$$
q_{-s}(\alpha)=\alpha \frac{p_{-s}(\alpha-s)}{\alpha -s}
$$
\emph{Proof}:
\begin{align*}
\tag*{\qed}\Gamma(s)\frac{q_{1-s}(\alpha)}{\alpha}=\int\limits_{0}^{+\infty} (-f(-t)e^t)^{s-1}e^{-\alpha t}dt=\Gamma(s)\frac{p_{1-s}(\alpha-s+1)}{\alpha -s+1}
\end{align*}
This way an approach through Laplace transform gives an intuitive explanation why
$$
q_n(\alpha)=\frac{\alpha p_n(\alpha+n)}{\alpha+n}
$$
\textbf{Remark 3.6}. One can write $p_n(\alpha)$ in the following manner:
$$
\frac{1}{(n-1)!}\frac{p_n(\alpha)}{\alpha}=\res_{t=0} \frac{e^{\alpha t}}{f^n(t)}dt=\frac{1}{2\pi i} \oint_{0} \frac{e^{\alpha z}}{f^n(z)} dz
$$
to make this analogy clear.\\
It should also be mentioned that using this formal approach one can write down the general sum by integer shifts of the argument of $p_{s}(\alpha)$.\\
~\\
\textbf{Proposition 3.7}. (\emph{formal}). Suppose $f$ is entire and real-valued, $\mathcal{R}$ stands for Ramanujan summation and $\Im$ stands for imaginary part. Then the following formal equality holds true:
$$
\sum_{n \in \mathbb{Z}}^{\mathcal{R}} \frac{p_{1-s}(n+x)}{n+x}=-\frac{2\pi}{\Gamma(s)}\Im \sum_{n=1}^\infty (f(2\pi i n))^{s-1} e^{2\pi i nx-\pi i s}
$$
\emph{Proof}: Consider Hurwitz formula for Riemann $\zeta$-function:
$$
\sum_{n \in \mathbb{Z}}^{\mathcal{R}} (n+x)^{-s}=\frac{(2\pi)^s}{\Gamma(s)}\sum_{n=1}^\infty n^{s-1} \cos\left(2\pi nx-\frac{\pi s}{2}\right)
$$
Consider an action on both sides of this equality by operator $\left( f\left(\frac{\partial}{\partial x}\right)/\frac{\partial}{\partial x}\right)^{s-1}$. Its action on $\cos(x)$ is known according to Euler's formula and Taylor formula, and its action on the powers $(n+x)^{-s}$ is known according to (3.5). \qed
\begin{center}
	\textbf{$\mathfrak{T}$-chains}
\end{center}
\textbf{Remark 3.7}. Suppose $f \in x+x^2 \mathbb{C}[[x]]$. Consider the action of inverse logarithmic derivative operator on $f$, i.e. $\mathfrak{T} f\coloneqq 1/\mathrm{dlog}(f)= \dfrac{f}{f'}$. Then $\mathfrak{T} f\in x+x^2 \mathbb{C}[[x]]$. Also the sine function $\sin(x)$ has the image $\tan(x)$ under this action.\\
~\\
\textbf{Proposition 3.8}: Suppose
$$
e^{\alpha (f(x))^{inv}}=\sum_{n=0}^\infty \frac{p_n(\alpha)x^n}{n!};~~~~~ e^{\alpha (\mathfrak{T} f(x))^{inv}}=\sum_{n=0}^\infty \frac{t_n(\alpha)x^n}{n!}
$$
Then the following equality holds true:
$$
\left(f'\left(\frac{\partial}{\partial \alpha}\right)\right)^n \frac{p_n(\alpha)}{\alpha}=\frac{t_n(\alpha)}{\alpha} \eqno (3.6)
$$
\emph{Proof}: Since $f \in x+x^2 \mathbb{C}[[x]]$, we have $f' \in 1+x \mathbb{C}[[x]]$. It then follows that the action of this operator is well-defined on the ring of polynomials. Also it does not change the degree of any polynomial and preserves the highest term coefficient. Then the result follows directly from the representations of $p_n(\alpha)$ and $t_n(\alpha)$ as residues.\\
(note: this equality formally holds even for fractional indexes, since we have representations as Laplace transforms)\qed\\
~\\
\textbf{Proposition 3.9}: For an operator $\mathfrak{T}: x+x^2 \mathbb{C}[[x]] \to x+x^2 \mathbb{C}[[x]]$ there exists a well-defined inverse operator $\mathfrak{T}^{-1}: x+x^2 \mathbb{C}[[x]] \to x+x^2 \mathbb{C}[[x]]$ and
$$
\mathfrak{T}^{-1} f(x)=x \exp\left(\int_{0}^{x}\left[\frac{1}{f(t)}-\frac{1}{t}\right]dt\right)
$$
\emph{Proof}: Suppose $g/g'=f$ and $g'(0)=1$.\\

1)
$$
\int_{0}^{x}\left[\frac{1}{f(t)}-\frac{1}{t}\right]dt \in x\mathbb{C}[[x]], \text{since}~ f(0)=0, f'(0)=1.
$$

2)
$$
\int_{0}^{x} \left[\frac{g'(t)}{g(t)}-\frac{1}{t}\right] dt=\ln\frac{g(t)}{t}\Biggr|_{0}^{x}=\ln \frac{g(x)}{x g'(0)}=\ln \frac{g(x)}{x}
$$
Now take the exponent of both parts of this equality to obtain the desired result.\qed\\
~\\
\textbf{Remark 3.8}. The action of inverse operator may also be rewritten as
$$
\mathfrak{T}^{-1} f(x)=f(x) \exp\left(\int_{0}^{x}\left[\frac{1-f'(t)}{f(t)}\right]dt\right)
$$

So what would happen if we take the functions connected by the powers of $\mathfrak{T}$? Do they have interesting properties? Consider chains of the following form:

$$
\begin{CD}
...@<\mathfrak{T}<<\dfrac{f}{f'} @<\mathfrak{T}<< f @<\mathfrak{T}<< f(x)\exp\left(\displaystyle\int_{0}^{x} \left[\frac{1-f'(t)}{f(t)}\right]dt\right)@<\mathfrak{T}<<...
\end{CD}
$$

In general case such chains may be highly nontrivial, although sometimes they surprisingly have very consistent structure. But before investigating the stable chains, consider the following figure. In our case the transform $f(x) \to f(px)e^{-Ax}/p$ changes polynomials of binomial type in understandable way. Also we should think of the chain for $f(Ax)/A$ like it is the same as the one for $f(x)$, since $\mathfrak{T}^k f(Ax)/A = ((\mathfrak{T}^k f) \circ (Ax))/A$. That means that for general function $f(x)$ it is enough to look at the properties of the deformation $\cdot~e^{-x}$ by the exponent $e^{-x}$, not necessary to look at all deformations by other $e^{Ax}$.\\
$$
\begin{CD}
    @.                           \omega(x)               @.                    \varphi(x)\\
    @.                           @AA{inv}A                                 @AA{inv}A\\
...@<\mathfrak{T}<<\dfrac{f(x)}{f'(x)} @<\mathfrak{T}<< f(x)@<\mathfrak{T}<< f(x)\exp\left(\displaystyle\int_{0}^{x}\frac{1-f'(t)}{f(t)}dt\right)\\
    @.                            @.                                                        @VV{~\cdot\phantom{.}e^{-x}}V\\
...@<\mathfrak{T}<<\dfrac{f(x)}{f'(x)-f(x)} @<\mathfrak{T}<< f(x)e^{-x}@<\mathfrak{T}<< f(x)\exp\left(\displaystyle\int_{0}^{x}\frac{e^t-f'(t)}{f(t)}dt\right)\\
    @.                           @VV{inv}V                                  @VV{inv}V\\
    @.                 \omega\left(\displaystyle\frac{x}{1+x}\right)               @.                    \gamma(x): @. (xe^{\gamma(x)})^{inv}=xe^{-\varphi(x)}\\
\end{CD}
$$
~\\
There are two trivial statements on this figure which connect two operators $(\cdot)^{inv}$ and $(\cdot)e^{-x}$. Let's prove them and an additional one before considering the periodic chains.\\
~\\
\textbf{Proposition 3.10}:
\begin{align*}
\tag*{(3.8)}&1)~ (\mathfrak{T}(f(x)e^{-x}))^{inv}=(\mathfrak{T}f(x))^{inv} \circ \left(\frac{x}{1+x}\right)\\
~\\
\tag*{(3.9)}&2)~ x e^{(f(x)e^{-x})^{inv}}=\left( x e^{-f^{inv}(x)} \right)^{inv}\\
~\\
\tag*{(3.10)}&3)~ \mathfrak{T}^2 (f(x)e^{-x})=(\mathfrak{T}^2 f(x))(1-\mathfrak{T}f(x))
\end{align*}
~\\
\emph{Proof}: Let, as on figure $\omega\coloneqq(f/f')^{inv}$, $\varphi= f^{inv}$, $\gamma(x)=(f(x)e^{-x})^{inv}$.\\
\begin{align*}
1)~(\mathfrak{T}(f(x)e^{-x}))^{inv}&=\left(\frac{f(x)}{f'(x)-f(x)}\right)^{inv}=\left(\frac{x}{1-x}\circ \mathfrak{T}f(x) \right)^{inv}=\\
\tag*{\qed}&=(\mathfrak{T}f(x))^{inv} \circ \left( \frac{x}{1+x}\right)=\omega\left( \frac{x}{1+x}\right)\\
~\\
2)~~\phantom{\mathfrak{T}((f}\tag*{\qed}(xe^{\gamma(x)})^{inv}&=(f(\gamma(x))e^{-\gamma(x)}e^{\gamma(x)})^{inv}=\gamma^{inv}\circ f^{inv}=(f(x)e^{-x})\circ \varphi=xe^{-\varphi(x)}\\
~\\
3)~\phantom{(\mathfrak{T}(}\mathfrak{T}^2 (f(x)e^{-x})&=\mathfrak{T} \frac{f}{f'-f}=\frac{f(f'-f)}{f'^2 -ff''}=\frac{ff'}{f'^2-ff''} \left(1-\frac{f}{f'}\right)=\\
\tag*{$\blacksquare$}&=\left[\frac{f}{f'}\middle/ \left(\frac{f}{f'}\right)' \right] (1-\mathfrak{T}f)=(\mathfrak{T}^2 f(x))(1-\mathfrak{T}f(x))
\end{align*}

\newpage

\hrule

\section{Deformation of trivial $\mathfrak{T}$-chains \protect\phantom{$q$}}

\hrule

~\\
~\\
As it was already mentioned, the chains may be very nontrivial in general case. But there are a few examples which are interesting, because their deformations are not too chaotic, and so the Taylor coefficients of the inverse functions of some links may be written down explicitly. Consider the main examples:\\

1-periodic chain and its deformation:
$$
\begin{CD}
@. ...@<\mathfrak{T}<<x @<\mathfrak{T}<< x @<\mathfrak{T}<< x @<\mathfrak{T}<<...\\
@.    @.                            @.                                                        @VV{~\cdot\phantom{.}e^{-x}}V\\
...@<\mathfrak{T}<<x-x^2@<\mathfrak{T}<<\dfrac{x}{1-x} @<\mathfrak{T}<< xe^{-x}@<\mathfrak{T}<< e^{\mathrm{Ei}(x)-\gamma}@<\mathfrak{T}<<...
\end{CD}
$$

2-periodic chain:
$$
\begin{CD}
@. ...@<\mathfrak{T}<<\dfrac{e^{px}-1}{p} @<\mathfrak{T}<< \dfrac{1-e^{-px}}{p} @<\mathfrak{T}<<\dfrac{e^{px}-1}{p} @<\mathfrak{T}<<...
\end{CD}
$$

And non-periodic but stable chain with hyperbolic $sine$ and $tangent$ functions, appearing with deformation $\cdot ~e^{-x}$  from the previous one in case $p=2$:
$$
\begin{CD}
...@<\mathfrak{T}<< \dfrac{\sinh(2px)}{2p} @<\mathfrak{T}<< \dfrac{\tanh(px)}{p} @<\mathfrak{T}<< \dfrac{\sinh(px)}{p} @<\mathfrak{T}<<...
\end{CD}
$$

So one could ask if there is any other example of periodic $\mathfrak{T}$-chain. An answer to this question is provided in section 5. For now consider the following functions:
\begin{align*}
&\Delta_p \coloneqq \tfrac{e^{px}-1}{p}&& (\Delta_0=x)\\
&y_p\coloneqq \Delta_p e^{-x}&& (\text{notice that}~ y_1=1-e^{-x};~ y_2=\sinh(x);~ y_0=xe^{-x})\\
&\gamma_p\coloneqq y_p^{inv}&& (\gamma_1=-\ln(1-x);~ \gamma_2=\ln(x+\sqrt{1+x^2});~ \gamma_0=-W(-x))\\
&\omega_p\coloneqq (\mathfrak{T}y_p)^{inv}&& (\text{``arctangent''})\\
&T_p \coloneqq \mathfrak{T}^{-1}y_p&& (T_1=e^x-1;~ T_0~ \text{was already considered in section 2 as}~ T)\\
&\psi_p \coloneqq {T_p}^{inv}
\end{align*}
Notice that, according to (3.10), we have $\mathfrak{T}^2 y_p=\Delta_p(1-\Delta_{-p})$ and consider the following figure:
$$
\begin{CD}
@.    @.                                      @.                    \dfrac{1}{p}\ln(1+px)\\
@.    @.                           @.                                @AA{inv}A\\
@....@<\mathfrak{T}<<\Delta_{-p} @<\mathfrak{T}<< \Delta_p@<\mathfrak{T}<< \Delta_{-p}\\
@.    @.                            @.                                                        @VV{~\cdot\phantom{.}e^{-x}}V\\
...@<\mathfrak{T}<<\Delta_p(1-\Delta_{-p}) @<\mathfrak{T}<<\dfrac{y_p}{y_p'} @<\mathfrak{T}<< y_p@<\mathfrak{T}<< T_p\\
@.    @.                           @V{inv}VV                                  @V{inv}VV                                             @V{inv}VV\\
@.    @.                 \omega_p              @.                               \gamma_p                  @.           \psi_p
\end{CD}
$$

\noindent The basic properties of considered functions are listed below.\\
~\\
\textbf{Proposition 4.1} (see~\cite{Askey},~\cite{Dif},~\cite{Gl}):
\begin{align*}
&\gamma_p(x)=\sum_{n=1}^\infty \frac{(px)^n}{n}\binom{\frac{n}{p}}{n}~~~ e^{\alpha\gamma_p(x)}=\sum_{n=0}^\infty \frac{\alpha(px)^n}{\alpha+n}\binom{\frac{\alpha +n}{p}}{n}
\end{align*}
\emph{Proof}: Follows from (3.3) \qed\\
~\\
\textbf{Proposition 4.2}:
$$
\omega_p(x)=\sum_{n=1}^\infty \frac{x^n}{pn}((p-1)^n-(-1)^n)=\frac{1}{p}\ln\left(\frac{1+x}{1+(1-p)x}\right)
$$
\emph{Proof}: Follows from (3.8) \qed\\
~\\
\textbf{Proposition 4.3}:
$$
xe^{\gamma_p(x)}=\left(x(1+px)^{-\frac{1}{p}}\right)^{inv}
$$
\emph{Proof}: Follows from (3.9) \qed\\
~\\
\textbf{Observation 4.1}: $y_p(x)$ is a solution of differential equation of the second order $y''_p=(p-2)y'_p+(p-1)y_p$ with characteristic polynomial $\chi(x)=(x+1)(x+1-p)$. It then follows that
$$
(\mathfrak{T} y_p)'=1+(2-p)(\mathfrak{T} y_p)+(1-p)(\mathfrak{T} y_p)^2~ \Rightarrow~ \omega_p'(x)=\frac{1}{x^2\chi(\frac{1}{x})}
$$
\textbf{Observation 4.2}:
$$
(y'_p+y_p)(y'_p+(1-p)y_p)^{p-1}=1
$$
\textbf{Observation 4.3}: For suitable $x$, we have
$$
\omega_p(x)=\int_{0}^\infty \frac{y_p(xt)}{t}e^{-t}dt
$$
\textbf{Observation 4.4}:
$$
\omega_p(x\gamma'_p)=\gamma_p~ \Rightarrow~ x\gamma'_p=\left(x(1+x)^{-\frac{1}{p}}(1+(1-p)x)^{\frac{1-p}{p}}\right)^{inv}
$$
\textbf{Observation 4.5}:
$$
\gamma'_p=\left(e^{(p-1)\gamma_p}-x\right)^{-1}~~ \text{or equivalently}~~ \ln\gamma'_p=(1-p)\gamma_p+\ln(1+x\gamma'_p)
$$
\textbf{Observation 4.6}:
\begin{align*}
y_p(A+B)&=y'_p(A)y_p(B)+y_p(A)y'_p(B)+(2-p)y_p(A)y_p(B)\\
y'_p(A+B)&=y'_p(A)y'_p(B)+(p-1)y_p(A)y_p(B)
\end{align*}
\textbf{Observation 4.7}. Notice the following invariants:
\begin{align*}
&y_{\frac{p}{p-1}}(x)=(1-p)y_p\left(\frac{x}{1-p}\right) &&\gamma_{\frac{p}{p-1}}(x)=(1-p)\gamma_p\left(\frac{x}{1-p}\right)\\
~\\
&y_{-p}(x)=y_p(x)e^{-px}
\end{align*}
And so the following holds true:
$$
y_{\frac{p}{p+1}}(x)=(1+p)y_{-p}\left(\frac{x}{1+p} \right)=(1+p)y_{p}\left(\frac{x}{1+p} \right)e^{\frac{px}{p+1}}
$$
It then follows that even if we don't know the general expansion of $\gamma_p$, we still could obtain the expansion of $\gamma_0$ using the transforms $p\to p/(p+1)\to p/(2p+1)...$, since on each step the corresponding polynomials of binomial type change in a predictable way, according to (3.3).\\
Unfortunately, such an approach does not help in studying the functions related to $\psi_p$.\\
~\\
\textbf{Proposition 4.4}:
\begin{align*}
&T_{\frac{p}{p-1}}(x)=(1-p)T_p\left(\frac{x}{1-p}\right) &&\psi_{\frac{p}{p-1}}(x)=(1-p)\psi_p\left(\frac{x}{1-p}\right)\\
~\\
&T_{-p}(x)=T_p(x)e^{p(e^x-1)}
\end{align*}
\emph{Proof}: Act with $\mathfrak{T}$ on both sides of the proposed equalities and notice that $1+p y_p e^x=e^{px}$. Now the result follows from the previous observation. \qed\\
~\\
\textbf{Corollary 4.1}:
\begin{align*}
T_{\frac{1}{n}}(nx)&=n\Delta_{1}e^{\Delta_{1}+\Delta_{2}+...+\Delta_{n-1}}\\
T_{\frac{2}{2n+1}}((2n+1)x)&=2(2n+1)\tanh\left(\frac{x}{2}\right)e^{2\Delta_{1}+2\Delta_{3}+...+2\Delta_{2n-1}}
\end{align*}
\textbf{Observation 4.7}:
$$
T_n(x)=(e^x-1)\prod_{k=1}^{n-1}\left(\frac{e^{x+\frac{2\pi i k}{n}}-1}{e^{\frac{2\pi i k}{n}}-1}\right)^{e^{-\frac{2\pi i k}{n}}}
$$
Although the expressions become rather complicated, it should be mentioned that there are some intersections of this structure in the domain of inverse functions for small $n$. As it was already mentioned, $\gamma_0(x)=-W(-x)$.\\
~\\
\textbf{Proposition 4.5}:
\begin{align*}
&\psi_1(x)=\ln(1+x) &&\psi_{-1}(x)=\ln(1+W(x))\\
&\psi_2(x)=2\mathrm{arcth}\left(\frac{x}{2}\right) &&\psi_{-2}(x)=\ln\left(1+\left(\frac{x}{1+\frac{x}{2}}e^{2x}\right)^{inv}\right)\\
&\psi_{\frac{1}{2}}(x)=2\ln\left(1+W\left(\frac{x}{2}\right)\right)
\end{align*}
\emph{Proof}: The result follows immediately from Proposition (4.4), Corollary (4.1) and Observation (4.7). \qed\\
~\\
\textbf{Remark 4.1}. We have already met the expression for $\psi_{-1}(x)$ in (1.3), but it should be considered as a coincidence, since there exists a canonical generalization of (1.3), provided below, which does not relate to $\psi_p$ in a direct way. We will return two this topic soon. For now consider the following two technical theorems, which generalize the construction, considered in section 2 of this paper.
\newpage
\noindent\underline{\textbf{Theorem 4.1}} (exact formula for $e^{\alpha \psi_p}$):
$$
\frac{e^{\alpha \psi_p(x)}-1}{\alpha}=\sum_{n=1}^\infty \frac{x^n}{n}\sum_{k=1}^{n}\frac{p^{k-1}\alpha^{n-k}}{(n-k)!}\sum_{m=0}^{k-1}\frac{(-n)^m}{m!}\sum_{\substack{\sum_{1}^{m}q_i=k-1\\ q_i>0}} \frac{B_{q_1}(\frac{1}{p})...B_{q_m}(\frac{1}{p})}{q_1 q_1!...q_m q_m!}
$$
where $B_n(x)$ are Bernoulli polynomials.\\
\emph{Proof}: See Appendix C.\\
~\\
\textbf{Remark 4.2}. Although this expression seems rather complicated, according to (3.2), polynomials $\nu_n^p(\alpha)=n![x^n]e^{\alpha \psi_p(x)}$ have the following nice property:
$$
\alpha~\frac{\nu_n^p(\alpha+p-1)-\nu_n^p(\alpha-1)}{p}=n\nu_n^p(\alpha)
$$
\underline{\textbf{Theorem 4.2}}:
\begin{align*}
&\left[ \frac{d^n}{dx^n} T_p^{\alpha}\right] T_p^{-\alpha}=\sum_{k=1}^{n} \sum_{m=k}^{n} (-1)^{m-k} \alpha^k e^{(k+p(m-k))x}\Delta_p^{-m} A_{k, m}^n(p)\\
&\left[ \frac{d^n}{dx^n} T_p^{\alpha}\right] T_p^{-\alpha}=\sum_{k=1}^{n} \sum_{m=k}^{n} (-1)^{n-k} \alpha^k e^{kx}\Delta_p^{-m} B_{k, m}^n(p)
\end{align*}

where $A_{k, m}^n(p), B_{k, m}^n(p)$ are expressions which satisfy the following relations:
\begin{align*}
&A_{k, m}^{n+1}(p)=(p(m-k)+k) A_{k, m}^n(p)+(m-1)A_{k, m-1}^n(p)+A_{k-1, m-1}^n(p)\\
&B_{k, m}^{n+1}(p)=(pm-k) B_{k, m}^n(p)+(m-1)B_{k, m-1}^n(p)+B_{k-1, m-1}^n(p)
\end{align*}
\emph{Proof}: By induction as in (2.2). \hfill\ensuremath{\blacksquare}\\
~\\
\underline{\textbf{Theorem 4.3}}:
\begin{center}
\boxed{\ln \gamma'_p(x)=\sum_{n=1}^\infty \frac{x^n}{n}\sum_{k=0}^{n}\binom{\frac{n}{p}}{k}p^k(1-p)^{n-k}}
\end{center}
~\\
The proof is divided in a few parts.\\
~\\
\textbf{Proposition 4.6}: Suppose $f^{inv}=\varphi$, $(f(x)e^{-x})^{inv}=\gamma(x)$ and $e^{\alpha \varphi(x)}=\sum_{0}^\infty p_n(\alpha)\frac{x^n}{n!}$. Then the following holds:
$$
-\ln(1-Axe^{\gamma(x)})=\sum_{n=1}^\infty \frac{x^n}{n} \sum_{k=0}^{n-1} \frac{p_k(n)}{k!}A^{n-k} \eqno (4.6)
$$
\emph{Proof}:
\begin{align*}
-\ln(1-Axe^{\gamma(x)})&=\sum_{k=1}^\infty \frac{A^kx^ke^{k\gamma(x)}}{k}=\\
&=\sum_{k=1}^\infty \frac{A^kx^k}{k}\sum_{m=0}^\infty \frac{kp_m(k+m)x^m}{(k+m)m!}=\\
&=\sum_{n=1}^\infty x^n \sum_{k=0}^{n-1}\frac{A^{n-k}}{n-k}\frac{(n-k)p_k(n-k+k)x^m}{(n-k+k)k!}=\\
\tag*{\qed}&=\sum_{n=1}^\infty \frac{x^n}{n} \sum_{k=0}^{n-1} \frac{p_k(n)}{k!}A^{n-k}
\end{align*}
\textbf{Proposition 4.7}:
$$
\ln \gamma'_p(x)=-\ln(e^{-\gamma_p(x)}-(1-p)x)  \eqno (4.7)
$$
\emph{Proof}: Follows from Observation (4.5). \qed\\
~\\
\underline{\emph{Proof of the theorem}}. According to (4.6) and (4.7), the following holds:
\begin{align*}
\ln \gamma'_p(x)&=\gamma_p(x)-\ln(1-(1-p)xe^{\gamma_p(x)})=\\
&=\sum_{n=1}^\infty \frac{p^nx^n}{n}\binom{\frac{n}{p}}{n}+\sum_{n=1}^\infty \frac{x^n}{n}\sum_{k=0}^{n-1}\binom{\frac{n}{p}}{k}p^k(1-p)^{n-k}=\\
\tag*{$\blacksquare$}&=\sum_{n=1}^\infty \frac{x^n}{n}\sum_{k=0}^{n}\binom{\frac{n}{p}}{k}p^k(1-p)^{n-k}
\end{align*}
For $p=0$ this formula results in (1.3). Since the radius of convergence of this expansion is equal to $\exp(\frac{1-p}{p}\ln(1-p))$ one can use the same approach as in section 1 and raise the question about the properties of the following series:
$$
M_p(s)=\sum_{n=1}^\infty \frac{(1-p)^{\frac{1-p}{p}n}}{n^s}\sum_{k=0}^{n}\binom{\frac{n}{p}}{k}p^k(1-p)^{n-k}
$$
(note: $M_1(s)=\zeta(s)$, $M_2(s)=2^{-s}\zeta(s)$)\\
~\\
\underline{\textbf{Theorem 4.4}}: (special case for the real pole of $\ln\gamma'_p$)
\begin{align*}
p \in [0; 1) &\Rightarrow \lim_{x \to 1^{-}} \ln\gamma'_p((1-p)^{\frac{1-p}{p}}x)+\frac{1}{2}\ln(1-x)=\frac{p-2}{2p}\ln(1-p)-\frac{1}{2}\ln 2\\
p=1 &\Rightarrow \lim_{x \to 1^{-}} \ln\gamma'_p((1-p)^{\frac{1-p}{p}}x)+1\ln(1-x)=0
\end{align*}
\emph{Proof}: See Appendix C.\\
~\\
\textbf{Remark 4.3}. It then follows from this theorem that for $p \in [0; 1)$ we have:
$$
\lim_{n\to\infty} (1-p)^{\frac{1-p}{p}n}\sum_{k=0}^{n}\binom{\frac{n}{p}}{k}p^k(1-p)^{n-k}=\frac{1}{2}
$$
\textbf{Remark 4.4}. As it was already shown, the deformation of the chain $...\leftarrow\Delta_{p}\leftarrow...$ results in relation $\mathfrak{T}^2 y_p=\Delta_{p}(1-\Delta_{-p})=y_py'_pe^{(2-p)x}$. According to (3.6) one can always write down the exact formula for the polynomials corresponding to the power $\mathfrak{T}^k y_p$ since the series for $f'(x)$ is known on each step, although it becomes more and more difficult to do.\\
~\\
\underline{\textbf{Theorem 4.5}}:
\begin{align*}
\frac{e^{\alpha(\mathfrak{T}^0 y_p)^{inv}}-1}{\alpha}&=\sum_{n=1}^\infty \frac{p^{n-1}x^n}{n}\binom{\frac{\alpha+n}{p}-1}{n-1}\\
\frac{e^{\alpha(\mathfrak{T}^1 y_p)^{inv}}-1}{\alpha}&=\sum_{n=1}^\infty \frac{x^n}{n}\sum_{k=0}^{n-1} \binom{\frac{\alpha}{p}-1}{k}\binom{n}{k+1}p^k(p-1)^{n-1-k}\\
\frac{e^{\alpha(\mathfrak{T}^2 y_p)^{inv}}-1}{\alpha}&=\sum_{n=1}^\infty \frac{x^n}{n}\sum_{k=0}^{n-1} \binom{\frac{\alpha}{p}+n-1}{k}\binom{2n-k-2}{n-1}p^k(1-p)^{n-1-k}\\
\frac{e^{\alpha(y_py'_p)^{inv}}-1}{\alpha}&=\sum_{n=1}^\infty \frac{x^n}{n}\sum_{k=0}^{n-1} \binom{\frac{\alpha+2n-p}{p}}{k}\binom{2n-k-2}{n-1}p^k(1-p)^{n-1-k}
\end{align*}
\emph{Proof}: See Appendix D.\\
~\\
Going back to the deformation of 1-periodic chain, one should look closely to the Taylor series of the inverses.\\
$$
\begin{CD}
\dfrac{x-x^2}{1-2x}@<\mathfrak{T}<<x-x^2@<\mathfrak{T}<<\dfrac{x}{1-x} @<\mathfrak{T}<< xe^{-x}\\
    @V{inv}VV                           @V{inv}VV                                  @V{inv}VV                                             @V{inv}VV\\
 \dfrac{1+2x-\sqrt{1+4x^2}}{2}   @.                 \dfrac{1-\sqrt{1-4x}}{2}             @.                               \dfrac{x}{1+x}               @. -W(-x)   
\end{CD}
$$
The Taylor series are:
$$
x+\sum_{n=1}^\infty \binom{2n}{n}\frac{(-1)^n x^{2n}}{2(2n-1)}; ~~\sum_{n=1}^\infty \binom{2n}{n}\frac{x^n}{2(2n-1)}; ~~\sum_{n=1}^\infty (-1)^{n-1}x^n; ~~  \sum_{n=1}^\infty \frac{n^{n-1}x^n}{n!}.
$$
~\\
\textbf{Remark}. Rational functions $\mathfrak{T}^k(x-x^2)=((e^{-x}\mathfrak{T})^{k+1} (e^x-1)) \circ (-\ln(1-x))$ have the property that the discriminants of their numerators form perfect squares. Also hypothetically for the exact exponent of 2 in factorization of discriminant the identity $\nu_2(D)=2^{2k-3}$ holds. This question is discussed in more detail on mathoverflow by user \href{https://mathoverflow.net/questions/296328/}{``Asymptotiac K''}. (https://mathoverflow.net/questions/296328/)
\begin{center}
\textbf{General comment on section 4}
\end{center}

When one wants to obtain ``natural'' generalization of some identity, different methods of proving this identity should be considered as different ways leading to possibly different generalizations (as for example the classical problem of interpreting natural logarithm as an integral of function $x^{-1}$ or as invertible power series at point 1). The method used in section 1 to obtain the expression for $\ln \gamma'_0$ was directly related to the similarity:
$$
n^n=\int_{0}^\infty t(t+n)^{n-1}e^{-t}dt~~ \sim~~ n!=\int_{0}^\infty t^n e^{-t}dt
$$
In general case an integral $\int p_n(t)e^{-t}dt$ for binomial polynomials $p_n$ is related to function $(1-\varphi(x))^{-1}$, and such an approach could only allow us to derive different expressions for $(\mathfrak{T}^k y_0)^{inv}$ in case of small $k$, not the general case $y_p$. As another example, to generalize the classical expansion of $\arcsin^2(x)$ one can consider the series related to \mbox{$\ln T_p( \ln(x(1+x^p)^{-1/p}))$}:
$$
\sum_{n=1}^\infty \frac{(-1)^{n-1}x^{(n-1)p}\Gamma\left(\frac{1}{p}\right)\Gamma(n)}{\Gamma\left(\frac{1}{p}+n\right)}
$$
Or use another appoach noticing that the following formal equality holds true:
$$
\sum_{n=1}^\infty \frac{x^n}{np^n \binom{\frac{n}{p}}{n}}=\int_{0}^\infty \frac{xe^{-t}}{1-xy_p(t)}dt
$$
Each of these series provides different proof for the expansion of $\arcsin^2(x)$ thus leading to the different types of generalizations. It should also be mentioned that surprisingly the function $T_0$ could be met in literature and it has direct relation to the Dickman function.

\newpage

\hrule

\section{Uniqueness of periodic chains}

\hrule

~\\
~\\
\underline{\textbf{Theorem 5.1}}: Suppose $k \in \mathbb{N}$, $f \in x+x^2 \mathbb{C}[[x]]$. Then the following holds true:
\begin{align*}
\mathfrak{T}^{2k}f&=f ~~\Rightarrow~~ \mathfrak{T}^{2}f=f~ \text{and there exists}~ p: f=\frac{e^{px}-1}{p};\\
\mathfrak{T}^{2k+1}f&=f ~~\Rightarrow~~ \mathfrak{T}f=f,~ \text{i.e.}~ f(x)=x.
\end{align*}
In other words the periodic chains are either 1- or 2-periodic.\\
(note.: $p$ may be equal to 0.)\\
~\\
It should be mentioned before providing the proof of this theorem that the first $n$ coefficients of $f$ determine the first $n$ coefficients of $\mathfrak{T}f$ and so they determine the first $n$ coefficients of $\mathfrak{T}^k f$ for any $k$.\\
~\\
\textbf{Proposition 5.1}: Consider all nonzero coefficients of $f$ of higher-order terms and take $x^n$ as the least of them. Suppose $\mathfrak{T}^k f=f$, and such a nontrivial term $x^n$ exists. Then $k$ is an even number.\\
~\\
\emph{Proof}: Suppose $f \approx x+ \alpha x^n$, $\alpha \neq 0$ . Then
$$
\mathfrak{T}f \approx \frac{x+\alpha x^n}{1+n\alpha x^{n-1}} \approx (x+\alpha x^n)(1-n\alpha x^{n-1}) \approx x+\alpha (1-n)x^n
$$
So it means, that
$$
\mathfrak{T}^k f \approx x+\alpha(1-n)^k x^n
$$
$f=\mathfrak{T}^k f$ and so $\alpha=(1-n)^k \alpha$. Now if $k$ is an odd number, then we have a contradiction, since after division by $\alpha$ LHS and RHS of this identity have different signs. \qed\\
~\\
\textbf{Corollary 5.1}: If the chain is periodic with odd period, then $f(x)=x$ and the part of theorem is proved.\\
\textbf{Corollary 5.2}: If the chain is periodic with even period and $f(x)\neq x$, then the coefficient $[x^2]f(x)$ is not equal to 0.\\
So for the function $\mathfrak{T}^{2k}f=f$ suppose that $f\approx x+\dfrac{x^2}{2!}+Ax^3+...$ without loss of generality.\\
~\\
\textbf{Proposition 5.2}. Suppose $n>1$.
\begin{align*}
\text{Suppose}~~~~~~~~~~f&\approx x+\frac{x^2}{2!}+\frac{x^3}{3!}+...+\frac{x^n}{n!}+\frac{\theta x^{n+1}}{(n+1)!}.\\
\text{Then}~~~~~~\mathfrak{T}^{2k}f &\approx x+\frac{x^2}{2!}+\frac{x^3}{3!}+...+\frac{x^n}{n!}+\frac{(1-n^{2k}(1-\theta)) x^{n+1}}{(n+1)!}.
\end{align*}
\emph{Proof}:
\begin{align*}
\mathfrak{T}f &\approx  \frac{x+\dfrac{x^2}{2!}+\dfrac{x^3}{3!}+...+\dfrac{x^n}{n!}+\dfrac{\theta x^{n+1}}{(n+1)!}}{1+x+\dfrac{x^2}{2!}+\dfrac{x^3}{3!}+...+\dfrac{x^{n-1}}{(n-1)!}+\dfrac{\theta x^{n}}{n!}}\approx\\
\end{align*}
\begin{align*}
&\approx \frac{x+\dfrac{x^2}{2!}+\dfrac{x^3}{3!}+...+\dfrac{x^n}{n!}+\dfrac{\theta x^{n+1}}{(n+1)!}}{e^x+\dfrac{\theta x^n}{n!}-\left(\dfrac{x^n}{n!}+\dfrac{x^{n+1}}{(n+1)!}+...\right)}\approx\\
&\approx e^{-x}\frac{x+\dfrac{x^2}{2!}+\dfrac{x^3}{3!}+...+\dfrac{x^n}{n!}+\dfrac{\theta x^{n+1}}{(n+1)!}}{1+\dfrac{x^n e^{-x}}{n!}\left((\theta-1)+\dfrac{x}{n+1}+...\right)}\approx\\
&\approx \left(x+\frac{x^2}{2!}+\frac{x^3}{3!}+...+\frac{x^n}{n!}+\frac{\theta x^{n+1}}{(n+1)!}\right)e^{-x}\left(1-\frac{(\theta-1)x^n}{n!}\right)\approx\\
&\approx  \left(x+\frac{x^2}{2!}+\frac{x^3}{3!}+...+\frac{x^n}{n!}\right)e^{-x}+\frac{(n+1-n\theta)x^{n+1}}{(n+1)!}\approx\\
\phantom{\mathfrak{T}f}&\approx 1-e^{-x}-\frac{x^{n+1}}{(n+1)!}+\frac{(n+1-n\theta)x^{n+1}}{(n+1)!}\approx\\
&\approx x-\frac{x^2}{2!}+\frac{x^3}{3!}+...+\frac{(-1)^{n-1}x^n}{n!}+\frac{(-1)^n x^{n+1}(1+(-1)^n n(1-\theta))}{(n+1)!}
\end{align*}
Using this operation $2k$ times we obtain the desired result. \qed\\
~\\
Now suppose $\mathfrak{T}^{2k}f=f$ and $\theta$ is the least coefficient of $x^{n+1}/(n+1)!$ which is not equal to $1$. Then $\theta=1-n^{2k}(1-\theta) \Rightarrow n^{2k}=1$ and that is impossible since $n>1$, $k>0$. It is left to say that the transform $f(x)\to f(px)/p$ commutes with $\mathfrak{T}$ to get the form $\frac{e^{px}-1}{p}$ for an arbitary $p$. \hfill\ensuremath{\blacksquare}\\

\begin{center}
\textbf{Conclusion}
\end{center}

The main goal of this study was to investigate the behaviour of functional inverses of the functions $f$ and $f/f'$. The solution of the periodic equation and the fact that it is the famous elementary function came as a surprise to me, since it is really hard to notice using only the definition of logarithmic derivative. The natural step is to investigate the ``eigenfunctions'' of operator $\mathfrak{T}$, such functions that $f(x)/f'(x)=f(px)/p$, although it goes far away from the initial problem, when one wants to understand the properties of associated polynomials of binomial type.

\newpage

\begin{center}
\textbf{Appendix A}
\end{center}

\hrule

~\\
~\\
\textbf{Proposition 2.6}.  $b_n \xrightarrow{n \to \infty} 0$.\\

\emph{Proof}: Suppose $T(x)=x\exp \left( \mathrm{Ei}(x)-\ln|x|-\gamma \right)=\psi^{inv}$.\\

1)
$$
\lim_{x \to -\infty} \frac{T(x)+e^{-\gamma}}{e^x/x}=-e^{-\gamma} \eqno (1)
$$
~\\

Then near the point $-\infty$ we have the following asymptotic relation:
\begin{align*}
T(x)&=-e^{\mathrm{Ei}(x)-\gamma} \Rightarrow\\
\tag*{\qed}&\Rightarrow e^{-\gamma}+T(x) \sim -e^{-\gamma}(e^{\mathrm{Ei}(x)}-1) \sim -e^{-\gamma}\mathrm{Ei}(x)\sim-e^{-\gamma}\frac{e^x}{x} 
\end{align*}

2)
$$
\lim_{t \to 0^{+}} t\psi'(-e^{-\gamma-t})=e^{\gamma} \eqno (2)
$$
~\\

After the change of variable $x \to \psi(x)$ the limit (1) may be rewritten as
\begin{align*}
\lim_{x \to (-e^{-\gamma})^{+}} (x+e^{-\gamma})x \psi'(x)&=-e^{-\gamma} \Rightarrow\\
&\Rightarrow \lim_{x \to 1^{-}} -e^{-\gamma}(x-1) \psi'(-e^{-\gamma} x)=1 \Rightarrow\\
\tag*{\qed}&\Rightarrow \lim_{t \to 0^{+}} te^{-\gamma}\psi'(-e^{-\gamma-t})=1
\end{align*}

3)
$$
\lim_{N \to \infty} \frac{b_1 +2 b_2+...+N b_N}{N}=1 \eqno (3)
$$
~\\

According to the identity $x e^{-\gamma}\psi'(-e^{-\gamma}x)=\sum_1^\infty nb_n x^n$ and (2), the following holds true:
$$
\text{As}~~ t \downarrow 0:~~ \sum_{n=1}^{\infty} nb_n e^{-nt} \sim \frac{e^t}{t} \sim \frac{1}{t}
$$

Since $b_n > 0$ it follows from Hardy-Littlewood tauberian theorem~\cite{HL} that
\begin{align*}
\tag*{\qed}\lim_{N \to \infty} \frac{b_1 +2 b_2+...+N b_N}{N}=1
\end{align*}

The result follows immediately. \hfill\ensuremath{\blacksquare}

\newpage

\begin{center}
\textbf{Appendix B}
\end{center}

\hrule

~\\
~\\
$n=1:$
\begin{gather*}
\tag*{$k=1$} 1
\end{gather*}
$n=2:$
\begin{gather*}
\tag*{$k=1$} 1 \qquad 1\\
\tag*{$k=2$}1
\end{gather*}
$n=3:$
\begin{gather*}
\tag*{$k=1$} 1 \qquad 2 \qquad 2\\
\tag*{$k=2$}3 \qquad 3\\
\tag*{$k=3$}1
\end{gather*}
$n=4:$
\begin{gather*}
\tag*{$k=1$} \phantom{1}1 \qquad \phantom{1}3 \qquad \phantom{1}6 \qquad \phantom{1}6 \\
\tag*{$k=2$}\phantom{1}7 \qquad 14 \qquad 11\\
\tag*{$k=3$}\phantom{1}6 \qquad \phantom{1}6\\
\tag*{$k=4$}\phantom{1}1
\end{gather*}
$n=5:$
\begin{gather*}
\tag*{$k=1$} \phantom{1}1 \qquad \phantom{1}4 \qquad 12 \qquad 24 \qquad 24\\
\tag*{$k=2$} 15 \qquad 45 \qquad 70 \qquad 50\\
\tag*{$k=3$} 25 \qquad 50 \qquad 35\\
\tag*{$k=4$} 10 \qquad 10\\
\tag*{$k=5$} \phantom{1}1
\end{gather*}
$n=6:$
\begin{gather*}
\tag*{$k=1$} \phantom{11}1 \qquad \phantom{11}5 \qquad \phantom{1}20 \qquad \phantom{1}60 \qquad 120 \qquad 120\\
\tag*{$k=2$} \phantom{1}31 \qquad 124 \qquad 287 \qquad 404 \qquad 274\\
\tag*{$k=3$} \phantom{1}90 \qquad 270 \qquad 375 \qquad 225\\
\tag*{$k=4$} \phantom{1}65 \qquad 130 \qquad \phantom{1}85\\
\tag*{$k=5$} \phantom{1}15 \qquad \phantom{1}15\\
\tag*{$k=6$} \phantom{11}1
\end{gather*}

\newpage

\begin{center}
\textbf{Appendix C}
\end{center}

\hrule

~\\
~\\
\underline{\textbf{Theorem 4.1}}.
$$
\frac{e^{\alpha \psi_p(x)}-1}{\alpha}=\sum_{n=1}^\infty \frac{x^n}{n}\sum_{k=1}^{n}\frac{p^{k-1}\alpha^{n-k}}{(n-k)!}\sum_{m=0}^{k-1}\frac{(-n)^m}{m!}\sum_{\substack{\sum_{1}^{m}q_i=k-1\\ q_i>0}} \frac{B_{q_1}(\frac{1}{p})...B_{q_m}(\frac{1}{p})}{q_1 q_1!...q_m q_m!}
$$
\emph{Proof}:
\begin{align*}
\frac{e^{\alpha \psi_p(x)}-1}{\alpha}&=\sum_{n=1}^\infty \frac{x^n}{n}\res_{t=0} \frac{e^{\alpha t}}{t^n}\exp\left(-n\int_{0}^{t}\left[\frac{pe^u}{e^{pu}-1}-\frac{1}{u}\right]du\right)dt\\
&=\sum_{n=1}^\infty \frac{x^n}{n}\sum_{k=1}^n \frac{\alpha^{n-k}}{(n-k)!}\res_{t=0} \frac{1}{t^k}\exp\left(-n\int_{0}^{t}\left[\frac{pe^u}{e^{pu}-1}-\frac{1}{u}\right]du\right)dt\\
&=\sum_{n=1}^\infty \frac{x^n}{n}\sum_{k=1}^n \frac{\alpha^{n-k}}{(n-k)!}\res_{t=0} \frac{1}{t^k} \sum_{m=0}^{k-1} \frac{(-n)^m}{m!}\left(\sum_{q=1}^\infty \frac{(pt)^q B_{q}(\frac{1}{p})}{qq!}\right)^mdt\\
&=\sum_{n=1}^\infty \frac{x^n}{n}\sum_{k=1}^{n}\frac{p^{k-1}\alpha^{n-k}}{(n-k)!}\sum_{m=0}^{k-1}\frac{(-n)^m}{m!}\sum_{\substack{\sum_{1}^{m}q_i=k-1\\ q_i>0}} \frac{B_{q_1}(\frac{1}{p})...B_{q_m}(\frac{1}{p})}{q_1 q_1!...q_m q_m!}\\
\tag*{$\blacksquare$}
\end{align*}
\underline{\textbf{Theorem 4.4}}. The case $p=1$ is obvious. Suppose $p \in [0,1)$. Denote $\pi_p\coloneqq(1-p)^{\frac{1-p}{p}}$. Notice that $\gamma_p(\pi_p)=-\frac{1}{p}\ln(1-p)$.
\begin{align*}
&\lim_{x \to 1^{-}} \ln\gamma'_p(\pi_p x)+\frac{1}{2}\ln(1-x)=\lim_{x \to \pi_p^{-}} \ln\gamma'_p(x)+\frac{1}{2}\ln(1-\pi_p^{-1}x)=\\
&=\lim_{-\frac{1}{p}\ln(1-p)^{-}} -\ln y'_p(x)+\frac{1}{2}\ln(1-\pi_p^{-1}y_p(x))=\\
&=\lim_{x\to 0^{-}} -\ln\left( \frac{(1-p)\pi_p}{p}(1-e^{px})\right)+\frac{1}{2}\ln\left(\frac{pe^x-e^{px}+1-p}{p}\right)=\\
&=\lim_{x \to 0^{-}} -\frac{1}{p}\ln(1-p)+\frac{1}{2} \ln \frac{p^2e^x-pe^{px}+p-p^2}{1+e^{2px}-2e^{px}}=\\
\tag*{$\blacksquare$}&=\frac{p-2}{2p}\ln(1-p)-\frac{1}{2}\ln 2
\end{align*}

\newpage

\begin{center}
\textbf{Appendix D}
\end{center}

\hrule

~\\
~\\
\underline{\textbf{Theorem 4.5}}. The fourth identity follows from the third one, according to Remark (4.4). The first identity is already known. It is left to prove the two remaining ones.\\

1)
\begin{align*}
\frac{e^{\alpha(\mathfrak{T}^1 y_p)^{inv}}-1}{\alpha}&=\sum_{n=1}^\infty \frac{x^n}{n} \res_{t=0} \frac{e^{\alpha t}(1+(p-1)e^{pt})^n}{(e^{pt}-1)^n}dt\\
&=\sum_{n=1}^\infty \frac{x^n}{np} \res_{t=0} \frac{e^{\frac{\alpha t}{p}}(1+(p-1)e^{t})^n}{(e^{t}-1)^n}dt\\
&=\sum_{n=1}^\infty \frac{x^n}{n} \res_{t=0} \frac{(1+pt)^{\frac{\alpha t}{p}}(1+(p-1)t)^n}{t^n}dt\\
\tag*{\qed}&=\sum_{n=1}^\infty \frac{x^n}{n}\sum_{k=0}^{n-1} \binom{\frac{\alpha}{p}-1}{k}\binom{n}{k+1}p^k(p-1)^{n-1-k}\\
\end{align*}

2)
\begin{align*}
\frac{e^{\alpha(\mathfrak{T}^2 y_p)^{inv}}-1}{\alpha}&=\sum_{n=1}^\infty \frac{x^n}{n}\res_{t=0} \frac{e^{\alpha t}}{\Delta_p^n (1-\Delta_{-p})^n}dt\\
&=\sum_{n=1}^\infty \frac{x^n}{np}\res_{t=0} \frac{e^{\frac{\alpha t}{p}}p^{2n}}{(e^t-1)^n (p-1+e^{-t})^n} dt\\
&=\sum_{n=1}^\infty \frac{x^n}{np} \res_{t=0} \frac{(1+t)^{\frac{\alpha}{p}+n-1}p^{2n}}{t^n(p+(p-1)t)^n}dt\\
&=\sum_{n=1}^\infty \frac{x^n}{n} \res_{t=0} \frac{(1+pt)^{\frac{\alpha}{p}+n-1} (1+(p-1)t)^{-n}}{t^n}dt\\
&=\sum_{n=1}^\infty \frac{x^n}{n}\sum_{k=0}^{n-1} \binom{\frac{\alpha}{p}+n-1}{k}\binom{2n-k-2}{n-1}p^k(1-p)^{n-1-k}\\
\tag*{$\blacksquare$}
\end{align*}

\end{document}